%% file: convex0.tex
%March 17, 2003
\input amstex
\documentstyle{amsppt}
\magnification=1200 \hsize=13.8cm \catcode`\@=11
\def\NoLogo{\let\logo@\empty}
\catcode`\@=\active \NoLogo

\def\heat{\lf(\frac{\p}{\p t}-\Delta\ri)}

\def\lf{\left}
\def\ri{\right}

\def\p{\partial}

\def\R{\Bbb R}

\def\vp{\varphi}

\def\abb{{\alpha\bar\beta}}

\def \D {\Delta}

\documentstyle{amsppt}
\magnification=1200 \hsize=13.8cm \vsize=19 cm

\leftheadtext{Lei Ni} \rightheadtext{Ricci flow and nonnegative
curvature} \topmatter
\title{Ricci flow and nonnegativity of curvature}\endtitle

\author{Lei Ni\footnotemark }\endauthor
\footnotetext"$^{1}$"{Research partially supported by NSF grant
DMS-0203023, USA.}
\address
Department of Mathematics, University of California, San Diego, La
Jolla, CA 92093
\endaddress
\email{ lni\@math.ucsd.edu}
\endemail

\affil { University of California, San Diego}
\endaffil

\date  April 2003\enddate

\abstract In this paper, we prove a general  maximum
principle for the time dependent Lichnerowicz heat equation on
symmetric tensors coupled with the Ricci flow on complete
Riemannian manifolds. As an application we construct complete
manifolds with bounded nonnegative sectional curvature of
dimension greater than or equal to four such that the Ricci flow does
not preserve the nonnegativity of the sectional curvature, even
though the nonnegativity of the sectional curvature was proved to
be preserved by Hamilton in dimension three. The example is the
first of this type. This fact is proved through a general
splitting theorem on the complete family of metrics with
nonnegative sectional curvature, deformed by the Ricci flow.
\endabstract

\endtopmatter

\document

\subheading{\S0 Introduction}\vskip .2cm

The Ricci flow has been proved to be an effective tool in the
study of the geometry and topology of manifolds. One of the good
properties of the Ricci flow is that it preserves the
`nonnegativity' of the curvature. In dimension three, Hamilton
\cite{ H1} proves that on compact manifolds the Ricc flow
preserves the nonnegativity of the Ricci curvature and the
sectional curvature. Using this property and the quantified
version, curvature pinching estimate,
 it was proved  in \cite {H1} that the normalized Ricci flow converges to
a Einstein metric if the initial metric admits positive Ricci
curvature. In particular, it implies that a simply-connected
 compact three-manifold is diffeomorphic to the three sphere if it admits a metric
with positive Ricci curvature. One can refer \cite{Ch} for an
updated survey  and \cite{P2} for 
some recent developement on the Ricci flow on three manifolds. Later in
\cite {H2} it was proved that the Ricci flow also preserves the
nonnegativity of the curvature operator in high dimension on
compact manifolds. In the K\"ahler case, Bando and Mok \cite{B, M}
proved that the flow also preserves the nonnegativity of the
holomorphic bisectional curvature. The Ricci flow on complete
manifold was initiated in \cite {Sh2}. In \cite{Sh3} Shi
generalized the above mentioned result of Bando and Mok to the
complete K\"ahler manifolds with bounded curvature. Interesting
applications were also obtained therein.

In this  paper, we shall study the topological consequences of the
assumption that Ricci flow preserves the nonnegativity of the
sectional curvature on complete Riemannian manifolds. The basic
method is to study the heat equation, time dependent, deformation
of the Busemann function via the optimal tensor maximum principle
proved in \cite{NT3}. The maximum principle of this type was first
proved by Hamilton for compact manifolds \cite{H2}. The proof of
\cite{H2} can be generalized to the complete noncompact manifolds
with bounded curvature with additional assumption that the tensor
 satisfying certain heat equation  is {\it uniformly bounded}
on the space-time. See for example \cite{NT2, Proposition 1.1}.
But in order to study the deformation of continuous functions, in
our case, the Busemann functions, one can not expect uniform
point-wise control on its Hessian since it is not even
differentiable. Therefore one needs a maximum principle assuming only
 integral bounds on the tensor considered. This is the main
technical difficulty. This difficulty was resolved in \cite{NT3}
and an optimal maximum principle was established there for the
time-independent heat equation. The tensor maximum principle
proved here in Theorem 2.1 is a time dependent version of the
 one  in \cite{NT3} (cf. Theorem 2.1 of \cite{NT3}).

By studying the deformation of the Busemann function, we shall
prove that  on a simply-connected complete Riemannian manifolds
with bounded nonnegative sectional curvature, if the Ricci flow
preserves the nonnegativity of the sectional curvature, then the
manifold splits as the product of a compact manifold with
nonnegative sectional curvature  with a complete manifold which is
diffeomorphic to the Euclidean space. As a corollary, we give
examples of complete Riemannian manifolds with {\it bounded}
nonnegative sectional curvature of dimension $\ge 4$ such that the
Ricci flow does not preserve the nonnegativity of the sectional
curvature. As far as we know, this is the first example of this
kind. Noticing that in dimension three the Ricci flow does
preserve the nonnegativity of the sectional curvature by \cite{H1}
on compact manifolds and complete manifolds with bounded
curvature. Another application of our approach is a classification
of complete manifolds with {\it bounded} nonnegative curvature
operator,  a result which has been previously established in \cite
{N} using different methods without assuming the boundedness of
the curvature (see also \cite{CY} for the compact case). The use
of the heat equation deformation of Busemann functions to study
the structure of complete manifolds was initiated in \cite{NT3}.
Therefore this paper can be viewed as a continuation of the
pervious work. The difference between this one and \cite{NT3} is
that we have to consider the heat equation with  metrics evolved
by the Ricci flow in order to have nice heat equation for the the
Hessian of its solution. Therefore we have to derive the heat
kernel estimate of Li-Yau type (cf. \cite{LY}) for the time
dependent heat equation. The estimate of this type was considered
before in \cite{Sa} by Saloff-Coste. However, the heat equation we
are considering does not belong to the classes considered in
\cite{Sa} (see Remark 1.1 for more details). Therefore we devote
the first section in establishing the heat kernel estimate as well
as the Harnack inequality for the time dependent heat equation,
following the approach of Grigoy'an in \cite{Gr1}. The result
itself might has its own interests.

\medskip

{\it Acknowledgement}. The author would like to thank Professors
Ben Chow, Peter Li and Jiaping Wang for their interests in this
work. Special thanks go to Professor Luen-Fai Tam for many helpful
suggestions. The paper is not possible without the previous
collaboration \cite{NT3} with him.

\newpage

\input convex1.tex
\input convex2.tex

\enddocument

%% file: convex1.tex
\input amstex
\documentstyle{amsppt}
\magnification=1200 \hsize=13.8cm \catcode`\@=11
\def\NoLogo{\let\logo@\empty}
\catcode`\@=\active \NoLogo
\def\heatt{\lf (\Delta-\frac{\p}{\p t}\ri)}

\def\heat{\lf(\frac{\p}{\p t}-\Delta\ri)}

\def\lf{\left}
\def\ri{\right}

\def\p{\partial}

\def\R{\Bbb R}

\def\vp{\varphi}

\def\abb{{\alpha\bar\beta}}

\def \D {\Delta}

\vsize=19.0 cm

\subheading{\S1 Time-dependent heat equation}

\vskip .2cm

Let $(M, g^0_{ij}(x))$ be a complete Riemannian manifold (of
dimension $n$) with bounded curvature tensor. We denote $k_0$ to
be the upper bound of $|R_{ijkl}|^2$, the curvature tensor of
$g^0$. By \cite{Sh2, Theorem 1.1, p. 224} we know that there
exists a constant $T(n, k_0)>0$ such that the Ricci flow
$$
\frac{\p}{\p t} g_{ij}(x,t)=-2R_{ij}(x,t) \tag 1.1
$$
has solution on $M\times [0, T]$. Moreover, there exists
$A'_m=A'_m(n, m, k_0)$ such that
$$
\|\nabla^m R_{ijkl}\|^2(x,t)\le \frac{A'_m}{t^m}. \tag 1.2
$$
In particular,
$$
\| R_{ijkl}\|(x,t)\le \sqrt{A_0} .\tag 1.3
$$
Moreover, $g_{ij}(x,t)$ has nonnegative curvature operator if the
initial metric $g_{ij}(x,0)$ has the nonnegative curvature
operator. We are going to study the initial value problem of the
heat equation
$$
\heatt v=0. \tag 1.4
$$
with initial vale $v(x,0)=u(x)$.  Here $\Delta v
=g^{ij}(x,t)v_{ij}$ with $v_{ij}$ denoting the Hessian of $v$.
Namely $\D$ is time-dependent. The following lemma is well-known
to experts. For example,  it was known and used in \cite{CH} by
Chow and Hamilton in their study of the linear trace Harnack
inequality for the Ricci flow.

\proclaim{Lemma 1.1}Let $v(x,t)$ be a solution to (1.4). Then the
complex Hessian $v_{ij}(x,t)$ satisfies
$$
\heat v_{ij}=2R_{ipjq}v_{pq}-R_{ip}v_{pj}-R_{pj}v_{ip}.   \tag 1.5
$$
\endproclaim
\demo{Proof} Direct calculation, using formulae on page 274 of
[H1], one has that
$$
(v_{ij})_t=(v_t)_{ij}+\lf(\nabla_iR_{jk}+\nabla_{j}R_{ik}
-\nabla_{k}R_{ij}\ri)v_k.
\tag 1.6
$$
On the other hand, the commutator calculation shows that
$$
v_{ijkk}=v_{kkij}+\lf(-\nabla_{s}R_{ij}+\nabla_i
R_{js}+\nabla_{j}R_{is}\ri)v_s+R_{is}v_{sj}+R_{js}v_{is}-2R_{isjk}v_{sk}.
\tag 1.7
$$
Now using (1.4) we have  $(v_t)_{ij}=v_{kkij}$. Then lemma follows
from (1.6) and (1.7).
\enddemo

\proclaim{Corollary 1.2} Denote $\eta$ be the symmetric tensor
$v_{ij}$. Denote $\|\eta\|^2$ the norm of $v_{ij}$ with respect to
$g_{ij}(x,t)$. Then   $\exp(-2\sqrt{A_0} t)\|\eta\|(x,t)$ is a
subsolution of (1.4).
\endproclaim
\demo{Proof} Direct calculation shows that
$$
\split \heatt \|\eta\|^2 &\ge
-4R_{ipjq}\eta_{pq}\eta_{ij}+4R_{ip}\eta_{pk}\eta_{ik}+ 2\|\nabla
\eta\|^2-4R_{ip}\eta_{pk}\eta_{ik}\\
&\ge 2\|\nabla \eta\|^2-4\sqrt{A_0}\|\eta\|^2.
\endsplit
$$
Here we have used Lemma 1.1. The claim of the corollary follows
easily.
\enddemo

In the following we collect some fundamental results on solution
(subsolutions) of (1.4). Our basic assumption is (1.3). For the
purpose of the later section we  also assume $T\le 1$ and
$g_{ij}(x,0)$ has nonnegative Ricci curvature. By (1.1) and (1.3)
we know that
$$
C(n, A_0)g_{ij}(x,0)\le g_{ij}(x,t)\le g_{ij}(x,0). \tag 1.8
$$
Since $g_{ij}(x,0)$ has nonnegative Ricci curvature, by (1.8), for
any $0\le t\le T$, we still have the following Neumann type
Poincar\'e inequality.

\proclaim{Lemma 1.2} Let $(M, g_{ij}(x,t))$ be a solution to the
Ricci flow such that the initial metric $g_{ij}(x,0)$ has
nonnegative Ricci curvature. For any domain $\Omega\subset B_0(y,
R)$ and any Lipschitz function $\vp$ on $\overline{\Omega}$
vanishes on $\p \Omega$
$$
\int_{\Omega}|\nabla \vp|^2\ge
\frac{b}{R^2}\lf(\frac{V_x(R)}{|\Omega|}\ri)^{\beta}\int_\Omega
\vp^2 \tag 1.9
$$
for some positive constants $\beta$, $b$ which only depends on $n$
and $A_0$. Here $|\nabla \vp|^2$ is calculated using
$g_{ij}(x,t)$, while $|\Omega|$ and $V_x(R)$ are calculated using
$g_{ij}(x,0)$.
\endproclaim

\demo{Proof} The lemma follows easily from Theorem 1.4 of
\cite{Gr1}. The point is that only the weak form
Neumann-Poincar\'e inequality and the volume doubling property are
needed in the proof of Theorem 1.4 of \cite{Gr1}. Since
$g_{ij}(x,0)$ has nonnegative Ricci curvature these two properties
hold for $(M, g_{ij}(x,0)$. On the other hand, the metric
$g_{ij}(x,t)$ is equivalent to $g_{ij}(x,0)$. Therefore these two
sufficient properties preserve.
\enddemo

 The next result is a mean value inequality. The proof is
 just a modification of the one for the time-independent heat equation case in
 \cite{Gr1}. Note that it is known from \cite{H1} that the scalar
 curvature ${\Cal R}(x,t)$ of $g_{ij}(x,t)$ is nonnegative, under the assumption
that $g_{ij}(x,0)$ has nonnegative Ricci, therefore scalar,
curvature.

\proclaim{Theorem 1.1} Let $(M,g(t))$ be as in Lemma 1.2. Let
$w(x,t)$ be a smooth function satisfying
$$
\heatt w\ge 0\tag 1.10
$$
on $\coprod_{\sqrt{t}}$ with $t\le T$, where
$\coprod_R=B_0(x,R)\times (0, R^2)$ and $B_{\tau}(x,\sqrt{t})$ is
the ball of radius $\sqrt{t}$ with respect to $g_{ij}(x,\tau)$.
 Then
$$
w^2_{+}(x,t)\le \frac{C(n, A_0, T)}{V_x(\sqrt{t})
t}\int_0^t\int_{B_0(x, \sqrt{t})}w^2_{+}(y,\tau)\, dyd\tau \tag
1.11
$$
Here $w_{+}$ is the positive part of $v$.
\endproclaim
\demo{Proof} Here we basically follow the argument of the proof of
Theorem 3.1 in \cite{Gr1}. The  key to the argument is the fact
that $g_{ij}(x,t)$ satisfying the Neumann-Poincar\'e inequality
(1.9) and the volume double property.  We  have these two
properties if we assume that the initial metric has nonnegative
Ricci curvature. To make the iteration argument work using Lemma
1.2 we need also to prove that the  Lemma 3.1 of \cite{Gr1} still
holds. In fact, for any $R\le \sqrt{t}$, let $\phi(x,t)$ be a
cut-off function on $B_0(x,R)$ such that $\phi(x,0)=0$. For
$\theta >0$, let $w_{\theta}=(w-\theta)_{+}$.
Multiplying
$w_\theta \phi^2$ on both side of (1.10) we have that
$$
\split \int_{\{w\ge \theta\}}w_t w_\theta \phi^2 &\le \int_{\{w\ge
\theta\}}(\D w)w_\theta \phi^2\\
&=-2\int_{\{w\ge \theta\}}<\nabla w_\theta, \nabla \phi>w_\theta
\phi-\int_{\{w\ge \theta\}}|\nabla w_\theta|^2\phi^2\\
& = -\int_{\{w\ge \theta\}}|\nabla
(w_{\theta}\phi)|^2+\int_M|\nabla \phi|^2w^2_\theta.
\endsplit \tag 1.12
$$
Integrating the time variable and noticing that
 $\phi\in C_0^{\infty}(B_0(x,R))$ we have that
$$
\int_0^t\int_{B_0(x,R)} w_{\theta}(w_{\theta})_\tau\phi^2 \le
-\int_0^t\int_{B_0(x,R)}|\nabla(w_{\theta}\phi)|^2+\int_0^t\int_{B_0(
x,R)}|\nabla\phi|^2w^2_{\theta}.
$$
The left hand side equals to
$$
\frac{1}{2}\int_0^t\int_{B_0(x,R)}(w^2_{\theta})_\tau\phi^2
=\frac{1}{2}\int_{B_0(x,R)}w^2_{\theta}\phi
 \vert ^{t}_{0}+\int_0^t\int_{B_0(x,R)}w^2_{\theta}\lf(-\phi_\tau\phi+\frac{1}{2}{\Cal
R}(y,\tau)\phi^2\ri).
$$
Combining the above two inequalities and using the fact ${\Cal
R}\ge 0$ we have that
$$
\int_{B_0(x,R)}w^2_{\theta}\phi^2(y, t)\, dy
+2\int_0^t\int_{B_0(x,R)}|\nabla(w_{\theta}\phi)|^2 \le
2\int_0^t\int_{B_0(x,R)}w^2_{\theta}\lf(|\nabla
\phi|^2+|\phi\phi_{\tau}|\ri).\tag 1.13
$$
Similarly, one can prove Lemma 3.2 of \cite{Gr1}, noticing that
Lemma 1.2 holds for metric $g_{ij}(x,t)$. Then the iteration
scheme in \cite{Gr1} can be applied to complete the proof of the
theorem.
\enddemo

For the Harnack inequality, let $v$ be a positive solution to
(1.4) on $\coprod_{8R}$ where $\coprod_R=B_0(x,R)\times (0, R^2)$.

\proclaim{Theorem 1.2} Let $(M, g_{ij}(x,t)$ be as in Lemma 1.2. Then
there exists a constant $\gamma =\gamma(n, A_0)>0$ such that
$$
v(x,64R^2)\ge \gamma \sup_{B_0(x, R)\times (3R^3, 4R^2)}v. \tag
1.14
$$
\endproclaim
\demo{Proof} The proof follows similarly as the proof of Theorem
4.1 in \cite{Gr1}. Since Lemma 4.2-4.4 in \cite{Gr1}
are robust enough to be adapted to current situation we only need
to establish the following result which corresponds to Lemma 4.1
of \cite{Gr1}. We can assume that $\sup_{\widetilde{\coprod}}v=1$.
\enddemo

\proclaim{Lemma 1.3} Let $v(x,t)$ be a positive solution to (1.4)
in $\coprod_{2R}$ and set
$$
H=\{(x,t)\in {\coprod} _R: v(x,t)>1\}, \quad \quad
\widetilde{\coprod}_R=B_0(x,R)\times (3R^2, 4R^2).
$$
Then for any $\delta>0$ there exists $\epsilon=\epsilon (\delta,
A_0, n)$ such that if
$$
|H|\ge \delta|{\coprod} _R|, \tag 1.15
$$
then
$$
\inf_{\widetilde{\coprod}_R}v\ge \epsilon.
$$
Here $|H|$ and $|\widetilde{\coprod}_R|$ are  measured with respect to
the 
metric $g_{ij}(x,0)$.
\endproclaim
\demo{Proof} We have similar situation as in the proof of Theorem
1.1. The argument follows closely as in \cite{Gr1}. Let $h=\log
(1/v)$. It is easy to have that $ \heat h=-|\nabla h|^2$. For a
cut-off function $\phi(x)$, we have  that
$$
\split \frac{\p}{\p t}\lf(\int_{B_0(x,R)}h_{+}\phi^2\ri)&
=\int_{B_0(x,R)}(h_{+})_t\phi^2-\int_{B_0(x,R)}h_{+}\phi^2{\Cal
R}\\
&\le \int_{B_0(x,R)}(h_{+})_t\phi^2 \\
&\le \int_{B_0(x,R)}(\D h_{+})\phi^2-|\nabla h_{+}|^2\phi^2\\
&\le -\frac{1}{2}\int_{B_0(x,R)}|\nabla
h_{+}|^2\phi^2+2\int_{B_0(x,R)}|\nabla \phi|^2.
\endsplit
$$
This is the (4.3) of \cite{Gr1}. The rest of the proof follows
verbatim as in the proof of \cite{Gr1, Lemma 4.1}.
\enddemo
One has the following immediate corollary of the above theorem.
\proclaim{Corollary 1.2} Let $v(x,t)$ be a weak positive solution
to (1.4) on $M\times [0, T]$. Then
$$
\frac{v(y,s)}{v(x,t)}\le
\exp\lf(C\lf(\frac{r^2(x,y)}{t-s}+\frac{t}{s}+1\ri)\ri). \tag 1.16
$$
Here $C=C(\gamma)>0$.
\endproclaim

\demo{Proof} This was proved, for example in \cite{Mo, page
110-112}.
\enddemo

\proclaim{Theorem 1.3} Let $(M, g_{ij}(t))$ be as the above. Let
$H(x,y, t)$ be the minimal positive heat kernel of the heat
equation (1.4). Then there exist positive constants $C_1$, $C_2$
 and $D$ such that
$$
C_1\frac{1}{V_0(x,\sqrt{t})}\exp\lf(-D\frac{r^2(x,y)}{t}\ri)\le
H(x,y, t)\le
C_2\frac{1}{V_0(x,\sqrt{t})}\exp\lf(-D\frac{r^2(x,y)}{t}\ri). \tag
1.17
$$
Here $V_0(x, a)$ and $r(x,y)$ denote the volume of $B_0(x,a)$ and
distance between $x$ and $y$, with respect to $g_{ij}(x,0)$,
respectively. $D>4$ is a absolute constant. $C_i=C_i(n,D)$.
\endproclaim
\demo{Proof}  Let $H(x,y,t)$ be the minimal positive heat kernel
of $\frac{\p}{\p t}-\D$. It is easy to see that for any $t>0$,
$$
\int_M H(x,y,t)\, dy_0 \le 1.
$$
Here $dy_0$ is the volume element with respect to the metric at
time $t=0$.  Fix a point $z\in M$ and let $u(x,t)=H(x,z,t)$. Then
there exists a point $y\in B_0(z,2\sqrt{t})$ such that
$$
u(y, 2t)\le \frac{1}{V_0(z,2\sqrt{t})}.
$$
Applying the Harnack and the volume doubling property we have that
$$
u(z,t)\le \frac{C(n)}{V_0(z, \sqrt{t})}. \tag 1.18
$$
Therefore we have the upper bound for $H(x,x,t)$. The upper bound
in (1.17) follows from a general result of Grigor'yan \cite{Gr2,
Theorem 1.1}. The lower bound can be obtained using the argument
in \cite{Gr1, page 73}. Let $\phi$ be a cut-off function such that
$\phi=1$ on $B_0(y,\frac{1}{2}\sqrt{t})$ and $\phi=0$ outside
$B_0(y,\sqrt{t})$. Now define
$$
w(x,s)=\int_M H(x,y,s)\phi(y)\, dy_0
$$
for $s\ge 0$ and $w(x,s)\equiv 1$ for $s\le 0$. Then $w(x,s)$ is a
solution to the heat equation on
$B_0(y,\frac{\sqrt{t}}{2})\times(-\infty, \infty)$. Here we have
extend the metric to be $g_{ij}(x,0)$ for $s\le 0$. Applying the
Harnack inequality (1.16) we have that
$$
\split 1&=u(y,0)\le C(n)u(y,\frac{\sqrt{t}}{2})\\
&=C(n)\int_M H(y,z, \frac{t}{2})\phi(z)\, dz_0\\
&\le C(n) \int_{B_0(y,\sqrt{t})}H(y,z,\frac{t}{2})\, dz_0\\
&\le C(n)\int_{B_0(y,\sqrt{t})}H(y,y,t) \, dz_0\\
&\le C(n)H(y,y,t)V_0(y,\sqrt{t}).
\endsplit
$$
This gives the lower bound for $H(x,x,t)$. The general form in
(1.17) is just another application of the Harnack inequality, or
Corollary 1.2.
\enddemo

\proclaim{Remark 1.1} In \cite{Sa}, the above Theorem 1.2 and Theorem
1.3  were proved
for the parabolic operator of type $\frac{\p}{\p t}-L$, with
$$
L f =m^{-1}\text{div}\lf( m {\Cal A}(\nabla f)\ri),
$$
where $m$ is a measure {\it independent} of $t$, ${\Cal A}$ is a
measurable section of $\text{End}\, (T_M)$ which is uniformly
equivalent to the identity. The time dependent Laplacian operator
can only expressed in the above form with time dependent  measure
$\sqrt{\det (g_{ij}(x,t))}dx_1\wedge\cdots\wedge dx_n$. Therefore
one can not just apply the results of \cite{Sa} directly. This
also is the reason that we need the scalar curvature of
$g_{ij}(x,t)$ is nonnegative to make the argument
 work. One could also prove the above theorems following the
 approach of Moser as in \cite{Sa}.
\endproclaim

%% file: convex2.tex
\input amstex
\documentstyle{amsppt}
\magnification=1200 \hsize=13.8cm \catcode`\@=11
\def\NoLogo{\let\logo@\empty}
\catcode`\@=\active \NoLogo
\def\heatt{\lf (\Delta-\frac{\p}{\p t}\ri)}

\def\heat{\lf(\frac{\p}{\p t}-\Delta\ri)}

\def\lf{\left}
\def\ri{\right}

\def\p{\partial}

\def\R{\Bbb R}
\def\Z{\Bbb Z}

\def\vp{\varphi}

\def\abb{{\alpha\bar\beta}}

\def \D {\Delta}

\documentstyle{amsppt}
\vsize=19.0 cm

\subheading{\S 2  A maximum principle for tensors and its
 applications}

\vskip .2cm

In this section  we shall prove a maximum principle for the
symmetric tensors satisfying (1.5) under the assumption that $(M,
g_{ij}(x,t))$ has bounded nonnegative sectional curvature. Since
the argument is very close to that  in \cite{NT3} we will be
sketchy here.

Let $\eta_{ij}$ be a symmetric tensor satisfying (1.5). The basic
assumption on $\eta$ is that there exists a constant $a>0$ such
that
$$
\int_M \|\eta\|(x,0)\exp\lf({-a r^2(x) }\ri)\, dx <\infty \tag 2.1
$$
and
$$
\liminf_{r\to\infty}\int_0^T\int_{B_0(o, r)}
\|\eta\|^2(x,t)\exp\lf({-ar^2(x)}\ri)\, dx\, dt <\infty. \tag 2.2
$$
Here $\|\eta\|(x,t)$ is the norm of $\eta_{ij}(x,t)$ with respect
to metrics $g_{ij}(x,t)$. But $B_0(0,r)$ is the ball with respect
to the initial metric $g_{ij}(x,0)$ and $r(x)$ is the distance
from $x$ to a fixed point $o\in M$ with respect to the initial
metric. Due to the fact that the maximum principle for the
heat equation does not hold on complete manifolds in general, one
needs some growth conditions on the solutions to make it true. The
condition (2.2) is optimal by comparing to the example given in
\cite{J, page 211-213}. The classical example there is a solution
to the heat equation on $\R\times [0,\times \infty)$, which has
zero initial data. The violation of the uniqueness implies the
failure of the maximum principle for the solutions. The example
has growth, as $|x|\to \infty$, just faster than $\exp(ar^2(x))$. The condition (2.1)
is needed to ensure that the equation (1.5) does have a solution
indeed. It is also in the  sharp form.

Before we state our result, let us first fix some notations. Let
$\varphi:[0,\infty)\to[0,1]$ be a smooth function so that
$\varphi\equiv1$ on $[0,1]$ and $\varphi\equiv0$ on $[2,\infty)$.
For any $x_0\in M$ and $R>0$, let $\varphi_{x_0,R}$ be the
function defined by
$$
\varphi_{x_0,R}(x)=\varphi\lf(\frac{r(x,x_0)}{R}\ri).
$$
Again, $r(x,y)$ denotes the distance function of the initial metric. 
Let $f_{x_0,R}$ be the solution of
$$
\heat f=-f
$$
  with initial value $\varphi_{x_0,R}$. Then $f_{x_0,R}$ is defined for all $t$ and is positive and bounded for $t>0$. In fact
$$
f_{x_0,R}(x,t)=e^{-t}\cdot \int_M H(x,y,t)\vp_{x_0,R}(y)dy_0.
$$

We shall establish the following maximum principle.

\proclaim{Theorem 2.1} Let $(M, g_{ij}(x,t))$ be a complete
noncompact Riemannian manifolds satisfying (1.1)--(1.3), with
nonnegative sectional curvature. Let $\eta (x,t)$ be a symmetric
tensor satisfying (1.5) on $M\times [0,T]$ with
$0<T<\frac{1}{40a}$ such that $||\eta||$ satisfies (2.1) and
(2.2). Suppose at $t=0$, $\eta_{ij}\ge -bg_{ij}(x,0)$ for some
constant $b\ge0$. Then there exists $0<T_0<T$ depending only on
$T$ and $a$ so that the following are true.
 \roster
\item"{(i)}" $\eta_{ij}(x,t)\ge -bg_{ij}(x,t)$ for all $(x,t)\in
M\times[0,T_0]$. \item"{(ii)}" For any $T_0>t'\ge0$, suppose that there
is a point $x'$ in $M^m$ and there exist constants $\nu>0$ and
$R>0$ such that the sum of the first $k$   eigenvalues
$\lambda_1,\dots,\lambda_k$ of $\eta_{ij}$ satisfies
$$
\lambda_1+\dots+\lambda_k\ge -kb+\nu k\varphi_{x',R}
$$
for all $x$ at time $t'$. Then for all $t>t'$ and for all $x\in
M$, the sum of the first $k$ eigenvalues of $\eta_\abb(x,t)$
satisfies
$$
\lambda_1+\dots+\lambda_k\ge  -kb+\nu kf_{x',R}(x,t-t').
$$
\endroster
\endproclaim
\demo{Proof} The complex version of Theorem 2.1 was proved in
\cite{NT3, Theorem 2.1}. The key step  of the argument is to
construct the barrier
$$
h(x,t)=\int_M H(x,y,t)\|\eta\| (y, 0)\, dy_0 \tag 2.3
$$
and $h_R(x,t)$ below  to control $\|\eta\|(x,t)$ on big annulus. It is easy
to see that $h(x,t)$ is a solution to (1.4). Using Lemma 1.4, the
assumption (2.2) and the maximum principle of \cite {NT1} we have
that
$$
\exp{(-2\sqrt{A_0}t)}\|\eta\|(x,t)\le h(x,t). \tag 2.4
$$
Let  $A_0(o, r_1, r_2)$ denote the annulus $B_0(o,r_2)\setminus
B_0(o,r_1)$. For any $R>0$, let $\sigma_R$ be a cut-off function
which is $1$ on $A_0(o,\frac{R}{4}, 4R)$ and $0$ outside $A_0(o,
\frac{R}{8}, 8R)$. We define
$$
h_R(x,t)=\int_M H(x,y,t)\sigma_R(y)||\eta||(y,0)dy_0.
$$
Then $h_R$ satisfies the heat equation with initial data
$\sigma_R||\eta||$. Since the proof of Lemma 2.3 of \cite{NT3}
only uses the heat kernel upper bound it remains to be true due to
Theorem 1.3. In particular, $h(x,t)\le \exp(2\sqrt{A_0}T)
(h_R(x, t)+\tau(R))$ on $A_0(o, \frac{R}{2}, 2R)$ with $\tau(R)\to 0$ 
as $R\to \infty$. And $h_R(x,t)\to 0$ as $R\to\infty$ on any compact
subset of $M$. Now using
$$
\exp{(2\sqrt{A_0}T)}\lf(h_R(x,t)+\tau(R)\ri)
$$
as the barrier the proof
of Theorem 2.1 follows verbatim as the corresponding result in
\cite{NT3}. Notice that the key inequality (2.14) in \cite{NT3},
still holds under the assumption that $g_{ij}(x,t)$ has
nonnegative sectional curvature (see also the proof of  Theorem
2.2 following).
\enddemo

The similar maximum principle for the scalar heat equations is
relatively easy to prove. They also require an assumption as
(2.2). The time dependent case was first proved in \cite{NT1}
following the original argument for the time-independent case in
\cite{L}. As an application we have the following approximation
result.

\proclaim{Theorem 2.2} Let $(M, g_{ij}(x,t))$ be as above. Let
$u(x)$ be a Lipschitz continuous convex function satisfying
$$
|u|(x)\le C\exp\lf(ar^2(x)\ri) \tag 2.5
$$
for some positive constants $C$ and $a$. Let $v(x,t)$ be the
solution to the time-dependent heat equation (1.4). There exists
$T_0>0$ depending only on $a$ and there exists $T_0>T_1>0$ such
that the following are true. \roster \item"{(i)}" For $0<t\le
T_0$, $v(\cdot,t)$ is a smooth convex function (with respect to
$g_{ij}(x,t)$). \item"{(ii)}" Let
$$
\Cal K(x,t)=\{w\in T^{1,0}_x(M)|\  v_{ij}(x,t) w^i=0,\text{\rm \
for all $j$}\}
$$
be the null space of $v_{ij}(x,t)$. Then for any $0<t<T_1$, $\Cal
K(x,t)$ is a  distribution on $M$. Moreover the distribution is
invariant in time as well as under the parallel translation.
\endroster
\endproclaim

In order to prove the above theorem we need the following
approximation result due to Greene-Wu \cite{GW3, Proposition 2.3}.

\proclaim{Lemma 2.1} Let $u$ be a convex function on $M$. Assume
that $u$ is Lipschitz with Lipschitz constant $1$. For any $b>0$,
there is a $C^\infty$ convex function $w$ such that \roster
\item"{(i)}" $|w(x)-w(y)|\le r(x,y)$; \item"{(ii)}" $|w-u|\le b$
on $M$; and \item"{(iii)}" $w_\abb\ge -bg_\abb$ on $M$.
\endroster
\endproclaim

\demo{Proof of Theorem 2.2} Once we have Lemma 2.1 and Theorem 2.1,
the proof follows as  the proof of Theorem
2.1, Corollary 2.1 and Theorem 3.1 of \cite{NT3}. The key fact
 is that under the
assumption $K_{ijij}\ge 0$,  for a choice of the orthogonal
frame such  that for the tensor $\eta_{ij}$ 
diagonalized at a fixed point $(x_0, t_0)$ with its 
eigenvalues $\lambda_i$ of $\eta$  ordered as  $\lambda_1\le
\lambda_2\le\cdots\le \lambda_n$,
$$
\split \sum_{i,j=1}^k \bigg[2R_{ipjq}&\eta_{pq}-R_{ip}\eta_{pj} -R_{pj}\eta_{i p }\bigg]g^{ij}\\
&=2\lf(\sum_{i=1}^k\sum_{p=1}^n R_{ipip}\lambda_p-\sum_{i=1}^k R_{ii}\lambda_i\ri)\\
&=2\lf(\sum_{i=1}^k\sum_{p=1}^n R_{ipip}\lambda_p-\sum_{i=1}^k \sum_{p=1}^n R_{ipip}\lambda_i\ri)\\
&=2\lf(\sum_{i=1}^k\sum_{p=k+1}^m\lambda_p R_{ipip}-\sum_{i=1}^k\sum_{p=k+1}^mR_{ipip}\lambda_i\ri)\\
&=2\lf(\sum_{i=1}^k\sum_{p=k+1}^m R_{ipip}(\lambda_p-\lambda_i)\ri)\\
&\ge 0.
\endsplit
$$
\enddemo

The following is the main result on the structure of solutions to
the Ricci flow preserving the nonnegativity of the sectional
curvature.

\proclaim{Theorem 2.3} Let $(M, g_{ij}(x,t))$ be solution to the
(1.1) satisfying (1.3) with nonnegative sectional curvature.
Assume also that $M$ is simply-connected. Then $M$ splits
isometrically as $ M=N\times M_1$, where $N$ is a compact manifold
with nonnegative sectional curvature. $M_1$ is diffeomorphic to
$\R^k$. For the restriction of metric $g_{ij}(x,t)$ on $M_1$
with $t>0$, there is a strictly convex exhaustion functions on
$M_1$. Moreover, the soul of $M_1$ is a point and the soul of $M$
is $N\times\{o\}$, where $o$ is the soul of $M_1$.
\endproclaim
\demo{Proof} Let ${\Cal B}$ be the Busemann function on $M$, with
respect to the initial metric $g_{ij}(x,0)$. As it was proved in
\cite{CG, GW2} that ${\Cal B}$ is a convex Lipschitz function with
constant $1$. Also it is an exhaustion function on $M$. In fact
${\Cal B}(x)\ge cr(x)$ when $r(x)$ is sufficient large, for some
$C>0$, where $r(x)$ is the distance function to a fixed point
$o\in M$. Let $v(x,t)$ be the solution of (1.4) with $v(x,0)={\Cal
B}(x)$. Under the assumption that $K_{ijij}\ge 0$ is preserved
under the Ricci flow (1.1),  we know that $v(x,t)$ is convex by
Theorem 2.2. 
Applying Theorem 2.2 again we know that the null space of $v_{ij}(x,t)$
is a parallel distribution on $M$. By the simply-connectedness of
$M$ and the De Rham's decomposition theorem we know that $M$
splits as $M=N'\times M'_1$, where on $M_1$, $(v_{ij}(x,t))>0$ and
$v_{ij}\equiv0$ on N. Since $v(x,t)$ is strictly convex and
exhaustive on $M'_1$, by Theorem 3 (a) of \cite{GW2} we know that
$M'_1$ is diffeomorphic to $\R^{k'}$, where $k'=\text{dim}(M'_1)$.
We claim that $N'$ is compact. Otherwise, $v$ is not constant
since it is exhaustive on $N'$ since $v$ is  an exhaustion
function on $M$ by Corollary 1.4 of \cite{NT3}. Using
$v_{ij}\equiv 0$ on $N'$, the gradient of $v$ is a parallel vector
field, which  gives the splitting of $N'$ as $N'=N''\times \R$,
such that $v$ is constant on $N''$. By the exhaustion of $v$ again
we conclude that $N''$ is compact. Also $v$ is a linear function
on the flat factor $\R$. But we already know that $v$ is
exhaustive, which implies that $v\to +\infty$ on both ends of
$\R$. This is a contradiction. This proves that $N'$ is compact.
Let $N=N'$ and $M=M_1'$ we have the splitting for $t>0$. It is
also clear that there exists strictly convex exhaustion function
on $M_1$. As for the splitting at $t=0$ we can obtain by the
limiting argument. First we have the isometric splitting
$M=N\times M_1$ as above for some fixed $t_1>0$. On the other
hand, by \cite{H2, Lemma 8.2} we know that the distribution given
by the null space of $v_{ij}$ is also invariant in time. Therefore, 
the splitting
$M=N\times M_1$ also holds for $0<t\le t_1$. Now just
taking limit as $t\to 0$ we have the metric splitting of
$g_{ij}(x,0)$ on $M$ as $N\times M_1$.

As a consequence of the fact that there exist strictly convex
exhaustion function on $M_1$, we know that the soul of $M_1$ is a
point. The reason is that first  the restriction of $v(x,t)$ to
its soul will be constant since the soul is a compact totally
geodesic submanifold. On the other hand it is strictly convex if
the soul, which is a totally geodesic submanifold, has positive
dimension. The contradiction implies that the soul of $M_1$ is a
point for $t>0$. For the case $t=0$ the result follows by the
homotopy consideration. Assume that the soul is not a point.
Denote the soul by ${\Cal S}(M_1)$. Then since ${\Cal S}(M_1)$ is
the homotopy retraction of $M_1$ we know that $H_{s}(M_1)=\Z$,
where $s=\text{dim}({\Cal S}(M_1))\ge 1$. On the other hand since
we already know that $M_1$ is diffeomorphic to $\R^k$. Thus
$H_s(M_1)=\{0\}$, which is a contradiction. Therefore we know that
the soul of $M_1$ with respect to the initial metric is also a
point. The claim that the soul of $M$ is just $N\times\{o\}$
follows from the following simple lemma. \enddemo

\proclaim{Lemma 2.2} Let $N$ be a compact Riemannain manifolds
with nonnegative sectional curvature. Let $M_1$ be a complete
noncompact Riemannian manifold with nonnegative sectional
curvature. Let $M=N\times M_1$. Then the soul of $M$, ${\Cal
S}(M)=N\times {\Cal S}(M_1)$, where ${\Cal S}(M_1)$ is the soul of
$M_1$.
\endproclaim
\proclaim{Remark 2.1} Combining with Theorem 5.2 of \cite{NT3},
this in particular implies that if the $M$ is a complete K\"ahler
manifolds with nonnegative sectional curvature, whose universal
cover does not contain the Euclidean factor, then the soul of $M$
is either a point or the compact factor which is a compact
Hermitian symmetric spaces. In particular, the result holds if the
Ricci curvature of $M$ is positive somewhere.
\endproclaim

\demo{Proof of Lemma 2.2} For any point $z\in M$ we write 
$z=(x,y)$ according to the product. First of all, it is easy to see
that $N\times {\Cal S}(M_1)$ is totally geodesic. It is also
totally convex since any geodesic $\gamma(s)$ on $M$ can be
written as $(\gamma_1(s), \gamma_2(s))$, where $\gamma_i(s)$ are
geodesics in the factor. Therefore, due to the fact ${\Cal
S}(M_1)$ is totally convex we know that $\gamma(s)$ lies inside
$N\times {\Cal S}(M_1)$ if its two end points do.

Let $\gamma(s)$ be any geodesic ray issued from $p\in M$. (Let
$p=(x_0, y_0)$ according to the product.) Since $N$ is compact we
have that for the projection $\gamma(s)=(\gamma_1(s),
\gamma_2(s))$, $\gamma_1(s)=p$ and $\gamma_2(s)$ is a ray in
$M_1$. Let ${\Cal B}^{\gamma}$ be the Busemann function with
respect to $\gamma$. We claim that ${\Cal B}^{\gamma}(x,y)={\Cal
B}^{\gamma_2}(y)$, where ${\Cal B}^{\gamma_2}(y)$ is the Busemann
function of $\gamma_2$ in $M_1$. Once we have the claim we
conclude that the level set of ${\Cal B}^{\gamma}$ is just
$N\times$ the level set of ${\Cal B}^{\gamma_2}$ in $M_1$ and the
half space $H^{\gamma}=\{z\in M\, |\,  {\Cal B}^{\gamma}(z)\le
0\}$, as proved in \cite{LT, Proposition 2.1}, $H^{\gamma}=N\times
H^{\gamma_2}$. Since this is true for any ray we have that
$C=\bigcap_{\gamma} H^{\gamma}$ is given by $N\times C_{M_1}$,
where $C_{M_1}$ denote the corresponding totally convex compact
subset in $M_1$. As in \cite{CG}, if the compact totally convex
subset $C$ has non-empty boundary we define $C^a=\{ p, d(p, \p C) \ge
a\}$. It is easy to see that $C_a=N\times C^a_{M_1}$. In
particular, this implies that the soul of $M$ is just $N\times
{\Cal S}(M_1)$ since the soul of $M$ is constructed by retracting
$C^a$ iteratively.

Now we verify the claim ${\Cal B}^{\gamma}(x,y)={\Cal
B}^{\gamma_2}(y)$.  By the definition we have that
$$
\split
 {\Cal B}^{\gamma}(x,y)&=\lim_{s\to \infty} s-d((x,y),
\gamma(s))\\
&=\lim_{s\to \infty}
s-\sqrt{d_N^2(x,x_0)+d^2_{M_1}(y,\gamma_2(s))}\\
&= \lim_{s\to \infty}(s-d_{M_1}(y,\gamma_2(s))+
\lf(\sqrt{d^2_{M_1}(y,\gamma_2(s))}-\sqrt{d_N^2(x,x_0)+d^2_{M_1}(y,\gamma_2(s))}\ri)\\
&=\lim_{s\to \infty} (s-d_{M_1}(y,\gamma_2(s))-\frac{d_N(x,x_0)}
{\sqrt{d^2_{M_1}(y,\gamma_2(s))}+\sqrt{d_N^2(x,x_0)+d^2_{M_1}(y,\gamma_2(s))}}\\
&=\lim_{s\to \infty} (s-d_{M_1}(y,\gamma_2(s))\\
&={\Cal B}^{\gamma_2}(y).
\endsplit
$$
This completes the proof of the lemma. 
\enddemo

Since the Ricci flow preserves the nonnegativity of the curvature
operator if the curvature is uniformly bounded (cf. \cite{H2}) we
have the following corollary on the structure of complete
simply-connected Riemannian manifolds with nonnegative curvature
operator.

\proclaim{Corollary 2.1} Let $M$ be a complete simply-connected
Riemannian manifold with bounded nonnegative curvature operator.
Then $M$ is a product of a compact Riemannian manifold with
nonnegative curvature operator with a complete noncompact manfold
which is diffeomorphic to $\R^k$. In the case of dimension three,
the same result holds if one only assumes that the  sectional
curvature is nonnegative.
\endproclaim

\proclaim{Remark 2.2} The compact factor in the above result has
been classified by Gallot and Meyer \cite {GaM} (also in \cite{CY}
by Chow and Yang) to be the product of compact symmetric spaces,
K\"ahler manifolds biholomorphic to the complex projective spaces
and the manifolds homeomorphic to spheres.

The above result was proved earlier in \cite{N} by Noronha without
assuming the boundedness of the curvature tensor. Our method here
has this restriction since we have to use the short time existence
result of Shi in \cite{Sh2} on  the Ricci flow. For dimension
three, in \cite{Sh1} the result was proved even for nonnegative
Ricci curvature case. However, it replies on the previous deep results
of Hamilton and Schoen-Yau. 

\endproclaim

As another application of Theorem 2.3 we give  examples of
complete Riemannian manifolds with nonnegative sectional curvature
on which the Ricci flow does not preserve the nonnegativity of the
sectional curvature. These manifolds can be constructed as
follows.  Let $G=SO(n+1)$ with the standard bi-invariant metric
and $H=SO(n)$ be its close subgroup. Then $H$ has action on $G$
(as translation) as well as its standard action on $P=\R^n$ (as
rotation). Let $M=G\times P/H$. Topologically $M$ is just the
tangent bundle over $S^n$ since $H\to G\to G/H=S^n$ is just the
corresponding principle bundle over $S^n$. These examples were
constructed in \cite{CG} to illustrate their structure theorem
therein. About these examples the following are known (cf.
\cite{CG}): The metric on $M$ has nonnegative sectional curvature
due to the fact that the metric is constructed as the base of
a Riemannian submersion; There is also another Riemannian
submersion from $T(S^n)$ to $S^n$ with fiber given by $\pi(g\times
P)$, where $\pi$ is the first submersion map from $G\times P$ to
$M$ (in general, there exists a Riemannian submersion from $M$ to
its soul according  to a result of Perelman \cite{P1}); The fiber
(which is given by $\pi(g\times P)$) of this submersion $\pi_*:
M\to S^n$ is totally geodesic; The fibers are not flat. Namely the
metric on each tangent space $T_p(S^n)$ is not the standard flat
metric; $M$ has the unique soul ${\Cal S}(M)=\pi(G\times 0)$ and
the metric on $M$ is not of product even locally.

\proclaim{Corollary 2.2} For the example manifolds above, the
Ricci flow does not preserve the nonnegativity of the sectional
curvature.
\endproclaim
\demo{Proof} First $M$ is simply-connected by the exact sequence
of the fibration $F\to M\to S^n$ with $F=\R^n$. Assume that the
Ricci flow preserves the nonnegativity of the sectional curvature.
By Theorem 2.3, we know that $M=N\times M_1$, where $M_1$ is
diffeomorphic to $\R^k$. This contradicts to the fact that the
metric on $M$ is not locally product (for most cases, it already 
contradicts to the fact that the tangent bundle $T(S^n)$ is
non-trivial topologically). In order to apply Theorem 2.3 we need to verify the
curvature of the initial metric is uniformly bounded. In the
following we focus on the case $M=SO(3)\times P/H$. The general
case follows from a similar consideration.

 As we know from \cite{CG, page 442 and CE, page 146-147},
 the metric is given such that $\pi: SO(3)\times \R^2\to T(S^2)$ is
a Riemannian submersion, where $SO(3)$ is the equipped with the
bi-invariant metric. Since the Riemannian submersion increases the
curvature, we know that  the metric constructed in this way has
nonnegative sectional curvature. The metric can also be described
using the second submersion $\pi_*$ from $T(S^2) \to S^2$ such
that for any point in the fiber if the tangent direction is
horizontal we use the metric from $S^2$ and  for the vertical
direction we use the metric given by
$$
dr^2+\frac{r^2}{1+r^2}d\theta^2.
$$
Here $(r,\theta)$ is the polar coordinates for $\R^2$. This
expression was claimed in \cite{CE, page 146}. For the sake of the
completeness we indicate the calculation here. Similar to the
situation considered in \cite{C, CGL} we can use $\frac{\p}{\p s}$
to denote the component of the Killing vector field of action
$SO(2)$ in $SO(3)$. The normalized Killing vector field is given
by
$$
W=\frac{1}{\sqrt{1+r^2}}\lf(\frac{\p}{\p \theta}+\frac{\p}{\p
s}\ri).
$$
Since
$$
\split
 {\Cal H}(\frac{\p}{\p \theta})& =\frac{\p}{\p
\theta}-<\frac{\p}{\p \theta}, W>W\\
&= \frac{\p}{\p \theta}-\frac{r^2}{1+r^2}W
\endsplit
$$
the metric on the base of $\frac{\p}{\p \theta}$ is given by
$$
\|\frac{\p}{\p \theta}\|^2_{M}=\|{\Cal H}(\frac{\p}{\p
\theta})\|^2=\frac{r^2}{1+r^2}.
$$
Here ${\Cal H}(\frac{\p}{\p \theta} )$ denotes the horizontal lift
(projection) of $\frac{\p}{\p \theta} $. This description make it
easy to verify that the curvature is uniformly bounded. In order
to calculate the curvature we need the formula of \cite{O'N} on
the submersion. The Corollary 1 of \cite{O'N, page 465} says that
$$
\split (a)\,  & K(P_{vw})=\hat K(P_{vw})
-\frac{<T_v v, T_w w>-\|T_v w\|^2}{\|v\wedge w\|^2} \\
(b)\, &K(P_{xv})\|x\|^2\|v\|^2=<(\nabla _xT)_v v,
x>+\|A_xv\|-\|T_v
x\|^2\\
(c)\, &K(P_{xy}) =K_{*}(P_{x_*y_*})-\frac{3\|A_x y\|^2}{\|x\wedge
y\|^2}, \text {  where  } x_*=\pi_*(x),
\endsplit \tag 2.6
$$
where $x,y$ are horizontal and $v,w$ are vertical. Here $A$ and
$T$ are the second fundamental form type tensor for the Riemannian
submersion $\pi_{*}:T(S^2)\to S^2$. $\hat K$ is the curvature of
the fiber and $K_{*}$ is the curvature of the base. Since the
fiber is totally geodesic, $T\equiv 0$,  we have the simplified
formula
$$
\split (a)\,  & K(P_{vw})=\hat K(P_{vw})\\
(b)\, &K(P_{xv})\|x\|^2\|v\|^2=\|A_xv\|^2\\
(c)\, &K(P_{xy}) =K_{*}(P_{x_*y_*}) -\frac{3\|A_x y\|^2}{\|x\wedge
y\|^2}, \text { where } x_*=\pi_*(x).
\endsplit \tag 2.7
$$
By (c) and the nonnegativity of $K(P_{xy})$ we have that
$K(P_{xy})$ is uniformly bounded. The curvature of the fiber can
be calculated directly. In fact in terms of the polar coordinates
on the fiber it is given by
$$
\frac{3}{(1+r^2)^2}.
$$
Therefore we have that $K(P_{vw})$ is also uniformly bounded. The
only thing need to be checked is the mixed curvature $K(P_{xv})$.
By the definition of $A$ we know that
$$
A_xv={\Cal H}\nabla_x V
$$
where ${\Cal H}$ is the horizontal projection and $V$ is any
arbitrary extension of $v$. For a unit horizonal vector $E$ we
have
$$
<A_x v, E>=-<v, \nabla_{x}E>.
$$
Therefore it is enough to show that the right hand side is
bounded. Since, by the first submersion consideration using the quotient, we
know that $K(P_{xy})$ is nonnegative. Therefore by (c) of (2.7),
$$
\|A_x y\|^2 \le \frac{1}{3}K_*(P_{x_* y_*})\|x\wedge y\|^2.
$$
This shows that $|<v, \nabla_{x}E>|$ is uniformly bounded. 

For the
sake of the completeness we also include a proof of the fact that
the fiber of $\pi '$ is totally geodesic since there is no written
proof in the literature. Recall that $M=SO(n+1)\times P/SO(n)$.
Here $SO(n)$ is viewed as close subgroup of $SO(n+1)$ by the
inclusion:
$$
A\to \lf(\matrix 1 & 0\cr 0& A\endmatrix\ri).
$$
We have the involution $\zeta$ which is given by
$$
\zeta=\lf(\matrix 1 &0\cr
                  0& -I \endmatrix\ri).
$$
$\zeta$ acts on $SO(n+1)$ by $A\to \zeta A\zeta$.  It is easy to
see that the fixed point set of $\zeta$ is $SO(n)$. Now we
consider the action of $\zeta$ on $SO(n+1)\times P$ as $(g, x)\to
(\zeta g\zeta, x)$. It is easy to see that this action is
commutative with the action of $SO(n)$ since for any $h\in SO(n)$
$\zeta h=h\zeta$. Therefore the action descends to $M$. It is easy
to see that the fixed point of this action is $\pi(e, P)$.  This
implies that the fiber $\pi(e,P)$ is totally geodesic since it is
the fixed point set of an isometry. The other fiber can be
verified similarly since for any point $p\in S^n$ there is also a
involution fixes that point.
\enddemo